\documentclass[english,10pt]{article}

\usepackage[latin1]{inputenc}
\usepackage{amsmath,amsthm,amssymb}
\usepackage{graphicx}
\usepackage{color}

\def\eq#1{(\ref{#1})}

\def\neweq#1{\begin{equation}\label{#1}}
\def\endeq{\end{equation}}

\def\ep{\varepsilon}
\def\o{\Omega}
\def\oo{\overline\Omega}

\def\RR{{\mathbb R} }

\def\di{\displaystyle}
\def\ri{\rightarrow}

\newtheorem{theorem}{Theorem}[section]
\newtheorem{prop}{Proposition}[section]
\newtheorem{lem}{Lemma}[section]
\newtheorem{corollary}{Corollary}[section]

\title{\sc A singular Gierer-Meinhardt system with different source terms}
\author{Marius GHERGU$^a$ and
Vicen\c tiu R\u ADULESCU$^{b,}$\thanks{Corresponding author}\\
\small $^a$ Institute of Mathematics ``Simion Stoilow" of the Romanian Academy,\\
\small 21, Calea Grivitei Street, 010702 Bucharest, Sector 1, Romania\\
\small E-mail: {\tt marius.ghergu@imar.ro}\\
\small $^b$ Department of Mathematics, University of Craiova,
200585 Craiova, Romania\\ \small E-mail: {\tt
vicentiu.radulescu@math.cnrs.fr}}

\date{}

\begin{document}
\baselineskip16pt 
\renewcommand{\theequation}{\arabic{section}.\arabic{equation}}
\catcode`@=11 \@addtoreset{equation}{section} \catcode`@=12

\maketitle

\begin{abstract}
We study the existence and nonexistence of classical solutions to
a general Gierer-Meinhardt system with Dirichlet boundary
condition. The main feature of this paper is that we are concerned
with a model in which both the activator and the inhibitor have
different sources given by general nonlinearities. Under some
additional hypotheses and in case of pure powers in
nonlinearities, regularity and uniqueness of the
solution in one dimension is also presented.\\
{\bf 2000 Mathematics Subject Classification}: 35J55, 35J65.\\
{\bf Key words}: Gierer-Meinhardt; Elliptic system;
Activator-inhibitor; Maximum principle.
\end{abstract}

\section{Introduction and the main results}

The systems of nonlinear equations of Gierer-Meinhardt type have
received a considerable attention in the last decade. These
problems arise in the study of biological pattern formation by
auto and cross catalysis being related to known biochemical
processes and cellular properties. The general model proposed by
Gierer and Meinhardt \cite{gmeinh,megi} may be written as
\neweq{gierr}
\left\{\begin{tabular}{ll} $\di u_t=d_1\Delta u-\alpha
u+c\rho\frac{u^p}{v^q}+\rho_0\rho$
\quad & $\mbox{\rm in}\ \Omega\times(0,T),$\\
$\di v_t=d_2\Delta v-\beta v+c'\rho' \frac{u^r}{v^s}$ \quad &
$\mbox{\rm in}\ \Omega\times(0,T),$
\end{tabular} \right.
\endeq
subject to Neumann boundary conditions. Here $\Omega\subset \RR^N$
$(N\geq 1)$ is a bounded domain, $u,v$ represent the
concentrations of the activator and inhibitor  with the source
distributions $\rho$ and $\rho'$ respectively. Also $d_1$, $d_2$
are diffusion coefficients with $d_1<<d_2$ and
$\alpha,\beta,c,c',\rho_0$ are positive constants. The exponents
$p,q,r,s\geq 0$ verify the relation $qr>(p-1)(s+1)>0$. The system
\eq{gierr} is of reaction-diffusion type and involves the
determination of an activator and an inhibitor concentration
field. In a biological context, the Gierer-Meinhardt system
\eq{gierr} has been used to model several phenomena arising in
morphogenesis and cellular differentiation.

The model presented by Gierer and Meinhardt \cite{gmeinh}
originates in the Turing's one \cite{turing} for morphogenesis in
the linear case and is based on the short range of activation and
on the long range of inhibition. Also the model introduced in
\cite{gmeinh} takes into account the classification between the
concentration of activators and inhibitors, on the one hand, and
the densities of their sources, on the other hand. A complete
description of entire dynamics of system \eq{gierr} is given in
the recent paper of Ni, Suzuki and Takagi \cite{nist}, where it is
shown that the dynamics of the system \eq{gierr} exhibit various
interesting behaviors such as periodic solutions, unbounded
oscillating global solutions, and finite time blow-up solutions.

Many recent works have been devoted to the study of the
steady-states solutions of \eq{gierr}, that is, solutions of the
stationary system
\neweq{gierrr}
\left\{\begin{tabular}{ll} $\di d_1\Delta u-\alpha
u+c\rho\frac{u^p}{v^q}+\rho_0\rho =0$
\quad & $\mbox{\rm in}\ \Omega,$\\
$\di d_2\Delta v-\beta v+c'\rho' \frac{u^r}{v^s}=0$ \quad &
$\mbox{\rm in}\ \Omega,$
\end{tabular} \right.
\endeq
\medskip
subject to Neumann boundary conditions. Such systems are difficult
to treat due to the lack of a variational structure or a priori
estimates. In this case it is more convenient to consider the {\it
shadow system} associated to \eq{gierrr}. More exactly, dividing
the second equation of \eq{gierrr} by $d_2$ and then letting
$d_2\rightarrow \infty$, we reduce the system \eq{gierrr} to a
single equation. The nonconstant solutions of such equation
present interior or boundary peaks or spikes, i.e., they exhibit a
{\it point  concentration} phenomenon. Among the great number of
works in this direction, we refer the reader to
\cite{nit2,nit1,niw, ww0, ww2, ww1} and the reference therein, as
well as to the survey paper of Ni \cite{ni}. For the study of
instability of solutions to \eq{gierrr}, we also mention here the
works of Miyamoto \cite{miy} and Yanagida \cite{yan1}.

In the case $\Omega=\RR^N$ ($N=1,2$) it has been shown in
\cite{dpino1,dpino2} that there exist ground state solutions of
\eq{giern} with single or multiple bumps in the activator which,
after a rescaling of $u$, are approaching a universal profile.

\medskip

Let $\Omega\subset \RR^N$ $(N\geq 1)$ be a bounded domain with
smooth boundary. In this paper we consider the stationary
Gierer-Meinhardt system for a wide class of nonlinearities subject
to homogeneous Dirichlet boundary conditions. More exactly, we are
concerned with the following elliptic system
$$
\left\{\begin{tabular}{ll} $\di \Delta u-\alpha
u+\frac{f(u)}{g(v)}+\rho(x)=0, \;u>0$
\quad & $\mbox{\rm in}\ \Omega,$\\
$\di \Delta v-\beta v+\frac{h(u)}{k(v)}=0, \;v>0$
\quad & $\mbox{\rm in}\ \Omega,$\\
$u=0,\; v=0$ \quad & $\mbox{\rm on}\ \partial\Omega,$
\end{tabular} \right.\eqno({\mathcal S})
$$
where $\alpha,\beta> 0$, $\rho\in C^{0,\gamma}(\Omega)$,
$(0<\gamma<1)$, $\rho\geq 0$, $\rho\not\equiv 0$ and $f,g,h,k\in
C^{0,\gamma}[0,\infty)$ are nonnegative  and nondecreasing
functions such that $g(0)=k(0)=0$. This last assumption on $g$ and
$k$, together with the Dirichlet conditions on $\partial\Omega$
makes the system singular at the boundary. Another difficulty is
due to the non-cooperative (i.e., non-quasimonotone) character of
our system.

We are mainly interested in the case where the activator and
inhibitor have different source terms, that is, the mappings
$t\longmapsto f(t)/h(t)$ and $t\longmapsto g(t)/k(t)$ are not
constant on $(0,\infty)$. Our study is motivated by some questions
addressed by Choi and McKenna \cite{choimck1, choimck2} or Kim
\cite{kimna, kimjmaa} concerning existence and nonexistence or
even uniqueness of the classical solutions for the model system
\neweq{giern}
\left\{\begin{tabular}{ll} $\di \Delta u-\alpha
u+\frac{u^p}{v^q}+\rho(x)=0$
\quad & $\mbox{\rm in}\ \Omega,$\\
$\di \Delta v-\beta v+\frac{u^r}{v^s}=0$
\quad & $\mbox{\rm in}\ \Omega,$\\
$u=0,\; v=0$ \quad & $\mbox{\rm on}\ \partial\Omega.$
\end{tabular} \right.
\endeq

In \cite{choimck1, kimna} it is assumed that the activator and
inhibitor have common sources and the approach rely on the
Schauder's fixed point theorem through a decouplization of the
system. More precisely, subtracting the two equations of
\eq{giern} we obtain in the case $p=r$ and $q=s$ a linear equation
in $w=u-v$. This  is suitable to obtain a priori estimates in
order to control the map whose fixed points are solutions of
\eq{giern}.

In Choi and McKenna \cite{choimck2} it is obtained the existence
of radially symmetric solutions of \eq{giern} in the case
$\Omega=(0,1)$ or $\Omega=B_1\subset\RR^2$ and $p=r>1$, $q=1$,
$s=0$. In \cite{choimck2} {\it a priori} bounds are obtained via
sharp estimates of the associated Green function.

In Section 2 we give a nonexistence result for classical solutions
to $({\mathcal S})$. To our best knowledge, there are no results
of this type in the literature. The main idea is to speculate the
asymptotic behavior of $v$ in the second equation of $({\mathcal
S})$. This will be then used in the first equation of the system
and by classical arguments  (see, e.g., \cite[Theorem 1.1]{grjde})
we obtain the desired nonexistence result. A special attention is
payed to the case of pure powers in nonlinearities. In this sense
we obtain some relations between the exponents $p,q,r$ and $s$ for
which the system \eq{giern} has no classical solutions.

In Section 3 we give an existence result for classical solutions
of $({\mathcal S})$ under the additional hypothesis $\beta \leq
\alpha$. In fact, this assumption is quite natural if we look at
the steady-state system \eq{gierrr}. We have only to divide the
first equation by $d_1$, the second one by $d_2$ and to take into
account the fact that $d_1<<d_2$. The existence in our case is
obtained without assuming any growth condition on $\rho$ near the
boundary since we are able to provide more general bounds for the
regularized system associated to $({\mathcal S})$. In particular,
we obtain that \eq{giern} has solutions provided that $r-p=s-q\geq
0$ and $q>p-1$.

The uniqueness of the solution  is a delicate matter. Actually,
there is only one result in the literature in this direction (see
\cite[Theorem 1]{choimck1}) and concerns the one dimensional case
of system \eq{giern} with $\rho \equiv 0$ and  $p=q=r=s=1$. Using
the same idea as in \cite{choimck1}, we are able to extend the
uniqueness of the solution to $({\mathcal S})$ in one dimension to
the following range of exponents: $0<q\leq p\leq 1$ and $r-p=s-q
\geq 0$. It is worth pointing out here that the uniqueness of the
solution for systems like $({\mathcal S})$ seems to be a
particular feature of the Dirichlet boundary conditions. As we can
see in the above mentioned works, in the case of Neumann boundary
conditions the Gierer-Meinhardt system does not have a unique
solution.

\section{A nonexistence result}
Several times in this paper we shall use the following result.  We
refer the reader to \cite[Lemma 2.1]{grjmaa} for a complete proof.
\begin{lem}\label{l1p}
Let $\Psi:{\Omega}\times(0,\infty)\rightarrow\RR$ be a H\"older
continuous function such that the mapping $\di(0,\infty)\ni
t\longmapsto\frac{\Psi(x,t)}{t}$ is strictly decreasing for each
$x\in\Omega.$ Assume that there exist  $v_1$, $v_2\in
C^2(\Omega)\cap C({\overline{\Omega}})$ such that

$(a)\quad \Delta v_1+\Psi(x,v_1)\leq 0\leq \Delta v_2+\Psi(x,v_2)$
in $\Omega;$

$(b)\quad v_1,v_2>0$ in $\Omega$ and $v_2\leq v_1$ on
$\partial\Omega;$

$(c)\quad\Delta v_1\in L^1(\Omega)\;\mbox{ or }\; \Delta v_2\in
L^1(\Omega).$

Then $v_2\leq v_1$ in $\Omega.$
\end{lem}

Another useful tool is the following result which is a direct
consequence of the maximum principle.

\begin{lem}\label{maxp}
Let $k\in C(0,\infty)$ be a positive nondecreasing function and
$a_1,a_2\in C(\o)$ with $0< a_2\leq a_1$ in $\o$. Assume that
there exist $\beta > 0$, $v_1,v_2\in C^2(\o)\cap C(\oo)$ such that
$v_1, v_2>0$ in $\o$, $v_1\geq v_2$ on $\partial\o$ and
$$\di \Delta v_1-\beta v_1+\frac{a_1(x)}{k(v_1)}\leq 0\leq
\Delta v_2-\beta v_2+\frac{a_2(x)}{k(v_2)}\quad\mbox{ in }\ \o.$$
Then $v_2\leq v_1$ in $\o$.
\end{lem}

Let $\Phi:[0,1)\ri [0,\infty)$ defined by
$$\di \Phi(t)=\int^t_0\frac{1}{\sqrt{2\int_\tau^1\frac{1}{k(\theta)}d\theta}}d\tau,\quad 0\leq t<1.$$

Set $a=\lim_{t\ri 1}\Phi(t)$ and consider $\Psi:[0,a)\ri [0,1)$
the inverse of $\Phi$. The main result of this section is the
following t nonexistence property.

\begin{theorem}
Assume that
\neweq{cint}
\di \int^a_0 \frac{tf(mt)}{g(M\Psi(t))}dt=+\infty,
\endeq
for all $0<m<1<M$. Then the system $({\mathcal S})$ has no
classical solutions.
\end{theorem}

{\it Proof.} Assume, by contradiction, that there exists a
classical solution $(u,v)$ of the system $({\mathcal S})$ and let
$\varphi_1$ be the normalized first eigenfunction of $-\Delta$ in
$H^1_0(\Omega)$. As it is well known, $\varphi_1 \in
C^2(\overline\Omega)$ and we can assume that $\varphi_1>0$ in
$\Omega$. Let  $\zeta$ denote the unique solution of the problem
\neweq{gier0}
\left\{\begin{tabular}{ll} $\di \Delta \zeta-\alpha
\zeta+\rho(x)=0$
\quad & $\mbox{\rm in}\ \Omega,$\\
$\di \zeta=0$ \quad & $\mbox{\rm on}\ \partial\Omega.$
\end{tabular} \right.
\endeq
By standard elliptic arguments and the classical maximum principle
we deduce that $\zeta\in C^2(\overline\Omega)$ and $\zeta>0$ in
$\Omega$.

In view of Hopf's maximum principle and taking into account the
regularity of the domain, there exist $c_1,c_2>0$ such that
\neweq{mzeta}
\di c_1 d(x)\leq \varphi_1,\zeta \leq c_2 d(x) \quad\mbox{ in }\
\Omega,
\endeq
where $d(x)=\mbox{dist}(x,\partial \Omega)$.

Since
$$
\left\{\begin{tabular}{ll} $\di \Delta
(u-\zeta)-\alpha(u-\zeta)\leq 0$
\quad & $\mbox{\rm in}\ \Omega,$\\
$\di u-\zeta = 0$\quad & $\mbox{\rm on}\ \partial\Omega,$\\
\end{tabular} \right.
$$
by the weak maximum principle  \cite[Corollary 3.2]{gt} we have
$u\geq \zeta$ in $\Omega$. Hence, by \eq{mzeta} it follows that
\neweq{dst}
\di u(x)\geq m d(x) \mbox{ in }\Omega,
\endeq
for some $m>0$ small enough. Set $C=\max_{x\in\overline\Omega}
h(u(x))>0$. Then $v$ satisfies
\neweq{inqqv}
\left\{\begin{tabular}{ll} $\Delta v-\beta v+\frac{C}{k(v)}\geq 0$
\quad & $\mbox{\rm in}\ \Omega,$\\
$v>0$ \quad & $\mbox{\rm in}\ \Omega,$\\
$v=0$ \quad & $\mbox{\rm on}\ \partial\Omega.$
\end{tabular} \right.
\endeq
Let $c>0$ be such that
\neweq{mrt}
\di c\varphi_1\leq \min \{a,d(x)\}\quad\mbox{ in }\Omega.
\endeq
We need the following auxiliary result.
\begin{lem}
There exists $M>1$ large enough such that $\overline v=M
\Psi(c\varphi_1)$ satisfies
\neweq{inqv}
\di \Delta \overline v-\beta\overline v+\frac{C}{k(\overline
v)}\leq 0\quad\mbox{ in }\Omega.
\endeq
\end{lem}
\noindent{\it Proof.} Since $\Phi(\Psi(t))=t$ for all $0\leq t<a$,
we get $\Psi(0)=0$ and $\Psi\in C^1(0,a)$ with
\neweq{der1psi}
\di \Psi'(t)=\sqrt{2\int_{\Psi(t)}^1\frac{1}{k(\tau)d\tau}}\,
\quad\mbox{ for all }0<t<a.
\endeq
This yields
\neweq{ps}
\left\{\begin{tabular}{ll} $\di-\Psi''(t)=\frac{1}{k(\Psi(t))}$
\quad & $\mbox{\rm for all}\ 0<t<a,$\\
$\Psi'(t), \Psi(t)>0$ \quad & $\mbox{\rm for all}\ 0<t<a,$\\
$\Psi(0)=0.$
\end{tabular} \right.
\endeq
By Hopf's maximum principle, there exist $\omega\Subset\Omega$ and
$\delta>0$ such that
\begin{equation}\label{omega}
\di |\nabla \varphi_1|>\delta \;\;\mbox{ in
}\;\Omega\setminus\omega \quad\mbox{ and
}\;\;\varphi_1>\delta\quad\mbox{ in }\;\omega.
\end{equation}

Fix $M>1$ large enough such that
\neweq{bm2}
M(c\delta)^2>C \quad\mbox{ and }\; Mc\lambda_1 \delta
\Psi'(c\|\varphi_1\|_\infty)>\frac{C}{\min_{x\in\omega}k(\Psi(c\varphi_1))}.
\endeq
We have
$$\di -\Delta \overline v=\frac{Mc^2}{k(\Psi(c\varphi_1))}|\nabla \varphi_1|^2
+Mc\lambda_1\varphi_1\Psi'(c\varphi_1)\quad\mbox{ in }\Omega.$$ By
\eq{bm2} we get
$$\begin{aligned}
\di -\Delta \overline v& \geq
Mc\lambda_1\varphi_1\Psi'(c\varphi_1)\geq
Mc\lambda_1\delta\Psi'(c\|\varphi_1\|_\infty)
\geq \frac{C}{k(\overline v)} &&\quad\mbox{ in }\omega,\\
\di -\Delta \overline v& \geq
\frac{Mc^2}{k(\Psi(c\varphi_1))}|\nabla \varphi_1|^2 \geq
\frac{C}{k(\Psi(c\varphi_1))}\geq \frac{C}{k(\overline v)}
&&\quad\mbox{ in }\Omega\setminus\omega.
\end{aligned}
$$

The last two inequalities imply that $\overline v$ satisfies
\eq{inqv}. This finishes the proof of the Lemma. \qed

\medskip

By virtue of Lemma \ref{maxp},  relations \eq{inqqv} and \eq{inqv}
yield $v\leq \overline v$ in $\Omega$. Using \eq{dst} we get
$$
\frac{f(u)}{g(v)}\geq \frac{f(md(x))}{g(M\Psi(c\varphi_1))} \quad
\mbox{\rm in}\ \Omega.
$$
Furthermore,  $u$ satisfies
\neweq{inqqv0}
\left\{\begin{tabular}{ll} $\di \Delta u-\alpha
u+\frac{f(md(x))}{g(M\Psi(c\varphi_1))}\leq 0$
\quad & $\mbox{\rm in}\ \Omega,$\\
$\di u>0$ \quad & $\mbox{\rm in}\ \Omega,$\\
$u=0$ \quad & $\mbox{\rm on}\ \partial\Omega.$
\end{tabular} \right.
\endeq
In order to avoid the singularities in \eq{inqqv0} near the
boundary, we consider the approximated problem
\neweq{inqqveps}
\left\{\begin{tabular}{ll} $\di \Delta w-\alpha
w+\frac{f(md(x))}{g(M\Psi(c\varphi_1))+\ep}= 0$
\quad & $\mbox{\rm in}\ \Omega,$\\
$\di w>0$
\quad & $\mbox{\rm in}\ \Omega,$\\
$w=0$ \quad & $\mbox{\rm on}\ \partial\Omega.$
\end{tabular} \right.
\endeq
Clearly $\overline w=u$ is a super-solution of \eq{inqqveps} while
$\underline w=0$ is a sub-solution. By classical results, the
problem \eq{inqqveps} has a unique solution $w_\ep\in
C^2(\overline\Omega)$ such that $w_\ep\leq u$ in $\Omega$.
Moreover, the maximum principle yields $w_\ep>0$ in $\Omega$.

In order to get a contradiction, we multiply by $\varphi_1$ in
\eq{inqqveps} and then we integrate over $\Omega$. We obtain
$$\di (\alpha+\lambda_1)\int_\Omega w_\ep \varphi_1 dx=\int_{\Omega} \varphi_1 \frac{f(md(x))}{g(M\Psi(c\varphi_1))+\ep} dx.$$
Since $w_\ep\leq u$ in $\Omega$ we have
$$\di (\alpha+\lambda_1)\int_\Omega u\varphi_1 dx\geq \int_{\omega}
\varphi_1 \frac{f(md(x))}{g(M\Psi(c\varphi_1))+\ep} dx\quad\mbox{
for all }\omega\Subset\Omega.$$ Let $\widetilde C=
(\alpha+\lambda_1)\int_\Omega u\varphi_1 dx$. Passing to the limit
in the above inequality we deduce
$$\di \int_{\omega} \varphi_1 \frac{f(md(x))}{g(M\Psi(c\varphi_1))} dx\leq \widetilde C<+\infty \quad\mbox{ for all }\omega\Subset\Omega.$$
Hence,
$$\di \int_{\Omega} \varphi_1 \frac{f(md(x))}{g(M\Psi(c\varphi_1))} dx\leq \widetilde C<+\infty.$$
Let now $\Omega_0=\{x\in\Omega: d(x)<a\}$. The above inequality
combined with \eq{mrt} produces
$$\di \int_{\Omega_0} d(x) \frac{f(md(x))}{g(M\Psi(d(x)))} dx<+\infty,$$
but this clearly contradicts \eq{cint}. Hence the system
$({\mathcal S})$ has no positive classical solutions. This
completes the proof of Theorem. \qed

\smallskip

If $k(t)=t^s$, $s>0$, condition \eq{cint} can be written more
explicitly by describing the asymptotic behavior of $\Psi$. We
have.

\begin{corollary}\label{cnon1}
Assume that $k(t)=t^s$, $s>0$, and one of the following conditions
hold
\begin{itemize}
\item[\rm (i)] $s>1$ and $\int^a_0 \frac{tf(mt)}{g(M
t^{2/(1+s)})}dt=+\infty,$ for all $0<m<1<M$; \item[\rm (ii)] $s=1$
and $\int^{\min\{a,1/2\}}_0 \frac{tf(mt)}{g(M t \sqrt{-\ln
t})}dt=+\infty,$ for all $0<m<1<M$; \item[\rm (iii)] $0<s<1$ and
$\int^a_0 \frac{tf(mt)}{g(M t)}dt=+\infty,$ for all $0<m<1<M$.
\end{itemize}

Then, the system $({\mathcal S})$ has no positive classical
solutions.
\end{corollary}

\noindent{\it Proof.} The main idea is to describe the asymptotic
behavior of $\Psi$ near the origin. Notice that in our case the
mapping $\Psi:[0,a)\ri [0,1)$ satisfies
\neweq{psi0}
\left\{\begin{tabular}{ll} $\di-\Psi''(t)=\Psi^{-s}(t)$
\quad & $\mbox{\rm for all}\ 0<t<a,$\\
$\Psi'(t), \Psi(t)>0$ \quad & $\mbox{\rm for all}\ 0<t<a,$\\
$\Psi(0)=0.$
\end{tabular} \right.
\endeq

\noindent (i) If $s>1$ then the mapping
$$\di (0,\infty)\ni t\longmapsto \left[\frac{(1+s)^2}{2(1-s)} \right]^{1/(1+s)}\cdot t^{2/(s+1)},$$
satisfies \eq{psi0}. Hence, there exist two positive constants
$c_1,c_2>0$ such that
$$\di c_1 t^{2/(s+1)}\leq \Psi(t)\leq c_2 t^{2/(s+1)}\quad\mbox{ for all }0<t<a.$$
Now, (i) follows directly from the above inequality.

\noindent (ii) Using the fact that $\Psi$ is concave, we deduce
that $\Psi(t)>t\Psi'(t),$ for all $0<t<a$. From \eq{psi0} it
follows that
$$\di -\Psi''(t)<\frac{1}{t\Psi'(t)}\quad\mbox{for all }\;0<t<a.$$
We multiply by $\Psi'$ in the last inequality and then we
integrate over $[t,b]$, $0<b<a$. We get
$$\di (\Psi')^{2}(t)-(\Psi')^{2}(b)\leq
2(\ln b-\ln t) \quad\mbox{for all }\;0<t\leq b<a.$$ Hence, there
exist $c_1>0$ and $\delta_1\in(0,b)$ such that $ \Psi'(t)\leq c_1
\sqrt{-\ln t}$ for all $0<t\leq \delta_1.$ Integrating over
$[0,t]$ we obtain
\neweq{HH3}
\di \Psi(t)\leq c_1\int_0^t \sqrt{-\ln \tau}d\tau=c_1t\sqrt{-\ln
t}+\frac{c_1}{2}\int_0^t \frac{1}{\sqrt{-\ln \tau}}d\tau
\quad\mbox{ for all } 0<t\leq \delta_1.
\endeq
Since the last integral in \eq{HH3} is finite, there exist $c_2>0$
and $\delta_2\in (0,\delta_1)$ such that
\neweq{HH4}
\di \Psi(t)\leq c_2 t\sqrt{-\ln t}\quad\mbox{for all }\;0<t\leq
\delta_2.
\endeq
From \eq{psi0} and \eq{HH4} we deduce
$-\Psi''(t)=\frac{1}{\Psi(t)}\geq \frac{1}{c_2}\cdot
\frac{1}{t\sqrt{-\ln t}}$ for all $0<t\leq\delta_2$.

An integration over $[t,\delta_2]$ in the last inequality yields
$$\di \Psi'(t)\geq \frac{2}{c_2}
\left( \sqrt{-\ln t}-\sqrt{-\delta_2}\right)\quad\mbox{for all
}\;0<t\leq \delta_2.$$ Therefore, there exist $c_3>0$ and
$\delta_3\in(0,\delta_2)$ such that $ \Psi'(t)\geq c_3 \sqrt{-\ln
t}$  for all $0<t\leq \delta_3.$ Proceeding in the same manner as
above, there exist $c_4>0$ and $\delta_4\in(0,\delta_3)$ such that
\neweq{HH5}
\di \Psi(t)\geq c_4 t\sqrt{-\ln t} \quad\mbox{for all}\;0<t\leq
\delta_4.
\endeq
From \eq{HH4} and \eq{HH5} we get
$$\di \di c_3 t\sqrt{-\ln t}\leq \Psi(t)\leq c_4 t\sqrt{-\ln t} \quad\mbox{for all}\;0<t\leq \delta_4.$$
Now, (ii) follows from the above estimates.

\noindent (iii) By \eq{der1psi} we have
$$\di \Psi'(t)=\sqrt{2\int_{\Psi(t)}^1 \tau^{-s}d\tau}=\sqrt{\frac{2}{1-s}(1-\Psi^{1-s}(t))}\, , \quad\mbox{ for all }0<t<a.$$
Hence $0<\Psi'(0)=\sqrt{2/(1-s)}<+\infty$ which implies $\Psi\in
C^1[0,a)$ and $c_1 t\leq \Psi(t)\leq c_2 t$ in $(0,a)$ for some
$c_1,c_2>0$. This proves (iii). \qed

\medskip

In the case of pure powers in the nonlinearities, we have the
following nonexistence result for \eq{giern}.

\begin{corollary}\label{cnon2}
Let $p,q,r,s>0$ be such that one of the following conditions hold
\begin{itemize}
\item[\rm (i)] $s>1$ and $2q\geq (s+1)(p+2)$; \item[\rm (ii)]
$s=1$ and $q>p+2$; \item[\rm (iii)] $0<s<1$ and $q\geq p+2$.
\end{itemize}

Then, the system \eq{giern} has no positive classical solutions.
\end{corollary}
\noindent{\it Proof.} The proofs of (i) and (iii) are simple
exercices of calculus. For (ii), by Corollary \ref{cnon2} we have
that \eq{giern} has no classical solutions provided $s=1$ and
\neweq{HH6}
\int_0^{1/2} t^{1+p-q}(-\ln t)^{-q/2} dt=+\infty.
\endeq
On the other hand, for $a,b\in\RR$ we have $\int_0^{1/2}
t^{a}(-\ln t)^{b} dt<+\infty$ if and only if $a>-1$ or $a=-1$ and
$b<-1$. Now condition \eq{HH6} reads $q>p+2$. This concludes the
proof. \qed

\section{Existence results}

For all $t_1,t_2>0$ we define
$$\di A(t_1,t_2)=\frac{f(t_1)}{h(t_1)}- \frac{g(t_2)}{k(t_2)}.$$
In this section we assume that $A$ fulfills

\medskip

\noindent $(A_1)\quad$ $A(t_1,t_2)\leq 0$ for all $t_1\geq t_2>0$.

\medskip

We also assume that

\medskip

\noindent $(A_2)\quad$ $k\in C^1(0,\infty)$ is nonnegative and
nondecreasing function such that
$\lim_{t\ri+\infty}\frac{K(t)}{h(t+c)}=+\infty$, for all $c>0$,
where $K(t)=\int_0^t k(\tau)d\tau$.

\medskip

Here are some examples of nonlinearities that fulfill $(A1)$ and
$(A2)$.
\begin{itemize}
\item[(i)] $f(t)=t^p,$ $g(t)=t^q$, $h(t)=t^r$, $k(t)=t^s$, $t\geq
0$, $p,q,r,s>0$, $r-p=s-q\geq 0$ and $p-q<1$;

\item[(ii)] $f(t)=\ln(1+t^p)$, $g(t)=e^{t^q}-1$, $h(t)=t^p$ and
$k(t)=t^q$, $t\geq 0$, $p,q>0$, $p-q<1$;
 \item[(iii)] $f(t)=\log(1+at)$,
$g(t)=\log(1+t)$, $h(t)=at$ and $k(t)=t$, $t\geq 0$, $a\geq 1$;
\end{itemize}

We give in what follows a general method to construct
nonlinearities $f,g,h,k$ that verify hypotheses $(A1)$ and $(A2)$.
Let $f,g,h,k:[0,\infty)\rightarrow [0,\infty)$ be nondecreasing
functions such that  $k$ and $h$ verify $(A2)$ and one of the
following assumptions hold:
\begin{itemize}
\item[$(a)$] $fk=gh$ and the mapping $(0,\infty)\ni t\longmapsto
f(t)/h(t)$ is nonincreasing; \item[$(b)$] there exists $m>0$ such
that $f(t)/h(t) \leq m\leq g(t)/k(t)$, for all $t>0$.
\end{itemize}
Then the mapping $A$ verifies $(A1)$.

For instance, the mappings given in example (i) satisfy the
condition $(a)$ while the mappings given in example (ii) verify
the condition $(b)$.

The first result of this section concerns the existence of
classical solutions for the general system $({\mathcal S})$.

\begin{theorem}\label{texi}
Assume that the hypotheses $(A1)-(A2)$ are fulfilled. Then the
system $({\mathcal S})$ has classical solutions.
\end{theorem}

The existence of a solution to $({\mathcal S})$ is obtained by
considering the regularized system
$$
\left\{\begin{tabular}{ll} $\di \Delta u-\alpha
u+\frac{f(u+\ep)}{g(v+\ep)}+\rho(x)=0$
\quad & $\mbox{\rm in}\ \Omega,$\\
$\di \Delta v-\beta v+\frac{h(u+\ep)}{k(v+\ep)}=0$
\quad & $\mbox{\rm in}\ \Omega,$\\
$u=0, \; v=0$ \quad & $\mbox{\rm on}\ \partial\Omega.$
\end{tabular} \right.\eqno{({\mathcal S})_\ep}
$$

\begin{lem}\label{lemep1} Let $u_\ep,v_\ep\in C^2(\Omega)\cap C(\overline\Omega)$
be a positive solution of $({\mathcal S})_\ep$. Then, there exists
$M>0$ which does not depend on $\ep$ such that
\neweq{conM}
\max\{\|u_\ep\|_\infty, \|v_\ep\|_\infty \}\leq M.
\endeq
\end{lem}

\noindent{\it Proof.} Let $w_\ep=u_\ep-v_\ep$ and
$\omega=\{x\in\Omega:w_\ep>0\}$. In order to prove the Lemma, it
suffices to provide an uniform upper bound for $v_\ep$ and
$w_\ep$. From $(S)_\ep$ we get
$$\begin{aligned}
\di \Delta w_\ep-\alpha w_\ep+\rho(x)&\di =(\alpha-\beta)v_\ep-
\frac{f(u_\ep+\ep)}{g(v_\ep+\ep)}+\frac{h(u_\ep+\ep)}{gk(v_\ep+\ep)}\\
&\di =(\alpha-\beta)v_\ep- \frac{h(u_\ep+\ep)}{g(v_\ep+\ep)}
A(u_\ep+\ep,v_\ep+\ep) \quad\mbox{ in }\Omega.
\end{aligned}$$
Let us notice that $A(u_\ep+\ep,v_\ep+\ep)\geq 0$ in $\omega$ and
$w_\ep=0$ on $\partial\omega$. This yields
$$
\di \Delta w_\ep-\alpha w_\ep \geq \rho(x) \quad\mbox{ in }\omega.
$$

Let $\zeta\in C^2(\oo)$ be the unique solution of \eq{gier0}. Then
$$
\left\{\begin{tabular}{ll} $\di \Delta
(w_\ep-\zeta)-\alpha(w_\ep-\zeta)\geq 0$
\quad & $\mbox{\rm in}\ \omega,$\\
$\di w_\ep-\zeta \leq 0$\quad & $\mbox{\rm on}\ \partial\omega.$\\
\end{tabular} \right.
$$
Furthermore, by the weak maximum principle  \cite[Corollary
3.2]{gt} we have $w_\ep\leq \zeta$ in $\omega$. Since $w_\ep\leq
0$ in $\Omega\setminus\omega$, it follows that
\neweq{mpalex1}
w_\ep\leq \zeta\quad\mbox{ in } \Omega.
\endeq

We  multiply by $k(v_\ep)$ in the second equation of $(S)_\ep$ and
we deduce that
\neweq{kdelta0}
\di k(v_\ep)\Delta v_\ep-\beta v_\ep
k(v_\ep)+\frac{k(v_\ep)}{k(v_\ep+\ep)}h(u_\ep+\ep)=0 \quad \mbox{
in } \Omega.
\endeq

But
\neweq{kdelta1}
\di k(v_\ep)\Delta v_\ep= \Delta K(v_\ep)-k'(v_\ep)|\nabla
v_\ep|^2 \quad \mbox{ in } \Omega.
\endeq
Since $k$ is nondecreasing, we have
\neweq{kdelta2}
\di K(v_\ep)=\int_0^{v_\ep}k(t)dt\leq v_\ep k(v_\ep)\quad\mbox{ in
} \Omega.
\endeq
Using now \eq{kdelta1} and \eq{kdelta2} in \eq{kdelta0} we deduce
$$
\di \Delta K(v_\ep)-k'(v_\ep)|\nabla v_\ep|^2- \beta
K(v_\ep)+\frac{k(v_\ep)}{k(v_\ep+\ep)}h(u_\ep+\ep)\geq 0 \quad
\mbox{ in } \Omega.
$$
Hence
\neweq{mpale}
\di \Delta K(v_\ep)-\beta K(v_\ep)+h(u_\ep+\ep)\geq 0 \quad \mbox{
in } \Omega.
\endeq

By \cite[Theorem 3.7]{gt}, there exists a positive constant $C>1$
depending only on $\Omega$ such that
$$\di \sup_{\overline\Omega}K(v_\ep)\leq  C\sup_{\overline\Omega}h(u_\ep+\ep)
\leq  C\sup_{\overline\Omega}h(v_\ep+\|\zeta\|_\infty+1).$$ Using
the assumption $(A2)$ we deduce that $(v_\ep)_\ep$ is uniformly
bounded, i.e., $\|v_\ep\|_\infty\leq m$ for some $m>0$ independent
on $\ep$. This yields $u_\ep=v_\ep+w_\ep$ $\leq
m+\|\zeta\|_\infty$ in $\Omega$ and the proof of Lemma
\ref{lemep1} is now complete. \qed

\begin{lem}\label{leme2}
For all $0<\ep<1$ there exists a solution $(u_\ep,v_\ep)\in
C^2(\oo)\times C^2(\oo)$ of the system $({\mathcal S})_\ep$.
\end{lem}

\noindent{\it Proof.} We use  topological degree arguments.
Consider the set
$$\di {\mathcal U}:=\left\{ (u,v)\in C^2(\oo)\times C^2(\oo): \begin{aligned}
\|u\|_\infty,\|v\|_\infty\leq M+1\\
u,v\geq 0 \mbox{ in }\o,
u\!\mid_{\partial\Omega}=v\!\mid_{\partial\Omega}=0
\end{aligned}
\right\},$$ where $M>0$ is the constant in \eq{conM}. Define
$$\Phi_t:{\mathcal U}\ri{\mathcal U},\quad \Phi_t(u,v)=(\Phi_t^1(u,v),\Phi_t^2(u,v)),$$
by
$$
\begin{aligned}
\Phi_t^1(u,v)&=u-t(-\Delta+\alpha)^{-1}\left(\frac{f(u+\ep)}{g(v+\ep)}+\rho\right),\\
\Phi_t^2(u,v)&=v-t(-\Delta+\beta)^{-1}\left(\frac{h(u+\ep)}{k(v+\ep)}\right).
\end{aligned}
$$

Using Lemma \ref{lemep1} we have $\Phi_t(u,v)\neq (0,0)$ on
$\partial{\mathcal U}$, for all $0\leq t\leq 1$. Therefore, by the
invariance of the topological degree at homotopy we have
$$\mbox{deg}\,(\Phi_1,{\mathcal U}, (0,0))= \mbox{deg}\,(\Phi_0,{\mathcal U},(0,0))=1.$$
Hence, there exists $(u,v)\in {\mathcal U}$ such that
$\Phi_1(u,v)=(0,0)$. This means that the system $({\mathcal
S})_\ep$ has at least one classical solution. \qed

\medskip

\noindent{\bf Proof of Theorem \ref{texi}.} Let $(u_\ep,v_\ep)\in
C^2(\oo)\times C^2(\oo)$ be a solution of $({\mathcal S})_\ep$.
  Then
$$
\left\{\begin{tabular}{ll} $\di \Delta
(u_\ep-\zeta)-\alpha(u_\ep-\zeta)\leq 0$
\quad & $\mbox{\rm in } \Omega,$\\
$\di u_\ep-\zeta = 0$\quad & $\mbox{\rm on }\partial\Omega,$\\
\end{tabular} \right.
$$
where $\zeta$ is the unique solution of \eq{gier0}. Hence
$\zeta\leq u_\ep$ in $\Omega$. By \eq{mpalex1} it follows that
\neweq{mpalex2}
\di w_\ep\leq \zeta\leq u_\ep\quad\mbox{ in } \Omega.
\endeq
Let $\xi\in C^2(\oo)$ be the unique positive solution of the
boundary value problem
\neweq{gz}
\left\{\begin{tabular}{ll} $\di \Delta \xi-\beta
\xi+\frac{h(\zeta)}{k(\xi+1)}=0$
\quad & $\mbox{\rm in}\ \Omega,$\\
$\di \xi=0$\quad & $\mbox{\rm on }\partial\Omega.$\\
\end{tabular} \right.
\endeq
In view of Lemma \ref{maxp} we have $\xi\leq v_\ep$ in $\o$, so
that, by Lemma \ref{lemep1}, the following estimates hold
\neweq{gz0}
\left\{\begin{tabular}{ll} $\di \zeta(x)\leq u_\ep(x)\leq M$
\quad & $\mbox{\rm in}\ \Omega,$\\
$\di \xi(x)\leq v_\ep(x)\leq M$
\quad & $\mbox{\rm in}\ \Omega.$\\
\end{tabular} \right.
\endeq
Now, standard H\"{o}lder and Schauder estimates can be employed in
order to deduce that $\{(u_\ep,v_\ep)\}_{0<\ep<1}$ converges (up
to a subsequence) in $C^2_{\rm loc}(\o)\times  C^2_{\rm loc}(\o)$
to $(u,v)\in C^2(\o)\times C^2(\o)$. It remains only to obtain an
upper bound near $\partial\o$ for $(u_\ep,v_\ep)$ which leads us
to the continuity up to the boundary of the solution $(u,v)$. This
will be done by combining standard arguments with the estimate
\eq{mpalex2}. First, by \eq{mpale} we have
\neweq{eqk}
\di \Delta K(v_\ep)+h(M+1)\geq 0 \quad \mbox{\rm in }\ \Omega.
\endeq
Fix $x_0\in\partial\o$. Since $\partial\o$ is smooth, there exist
$y\in\RR^N\setminus \o$ and $R>0$ such that $\oo\cap
\overline{B}(y,R)=\partial\o\cap \overline{B}(y,R)=\{x_0\}$.

Let $\delta(x)=|x-y|-R$ and $\o_0=\{x\in\o:4(N-1)\delta(x)<R\}$.

Consider $\psi\in C^2(0,\infty)$ such that $\psi'>0$ and
$\psi''<0$ on $(0,\infty)$ and set $\phi(x)=\psi(\delta(x))$,
$x\in\o_0$. Then
$$\begin{aligned}
\di \Delta \phi(x)&=\di \psi'(\delta(x))\Delta \delta(x)+\psi''(\delta(x))|\nabla \delta(x)|^2\\
&\di =\frac{N-1}{|x-y|} \psi'(\delta(x))+\psi''(\delta(x))\\
&\di \leq \frac{N-1}{R}
\psi'(\delta(x))+\psi''(\delta(x))\quad\mbox{ in }\o_0.
\end{aligned}$$
Let us choose now $\psi(t)=C\sqrt{t},$ $t>0$, where $C>0$.
Therefore
$$\di \Delta \phi(x)\leq \frac{C}{4}\delta^{-3/2}(x) \left[\frac{2(N-1)\delta(x)}{R}-1\right]
\leq -\frac{C}{8}\delta^{-3/2}(x)<0\quad\mbox{ in } \o.$$ We
choose $C>0$ large enough such that
\neweq{phi0}
\di \Delta \phi\leq -h(M+1)\quad\mbox{ in } \o_0
\endeq
and
\neweq{phi1}
\phi\mid_{\partial\o_0\setminus\partial \o}>K(M)\geq
\sup_{\overline\o_0} K(u_\ep).
\endeq
Furthermore, by \eq{eqk}, \eq{phi0} and \eq{phi1} we obtain
$$
\left\{\begin{tabular}{ll} $\di \Delta (\phi-K(v_\ep))\leq 0$
\quad & $\mbox{\rm in}\ \Omega_0,$\\
$\di \phi-K(v_\ep)\geq 0$
\quad & $\mbox{\rm on}\ \partial\Omega_0.$\\
\end{tabular} \right.
$$
This implies $\phi(x)\geq K(v_\ep)$ in $\o_0$, that is,
$$\di 0\leq v_\ep\leq K^{-1}(\phi(x))\quad\mbox{ in } \o_0.$$
Passing to the limit with $\ep\ri 0$ in the last inequality we
have $\di 0\leq v\leq K^{-1}(\phi(x))$ in  $\o_0.$ Hence
$$0\leq \lim_{x\ri x_0}v(x)\leq \lim_{x\ri x_0}K(\phi(x))=0.$$
Since $x_0\in\partial\o$ was arbitrary choosen, it follows that
$v\in C(\oo)$. Using the fact that $u_\ep=w_\ep+v_\ep\leq
\zeta+v_\ep$ in $\Omega$, in the same manner we conclude $u\in
C(\oo)$. This finishes the proof of Theorem \ref{texi}. \qed

\smallskip
The next result concerns the following system
\neweq{p}
\left\{\begin{tabular}{ll} $\di \Delta u-\alpha
u+\frac{u^p}{v^q}+\rho(x)=0$
\quad & $\mbox{\rm in}\ \Omega,$\\
$\di \Delta v-\beta v+\frac{u^{p+\sigma}}{v^{q+\sigma}}=0$
\quad & $\mbox{\rm in}\ \Omega,$\\
$u=v=0$ \quad & $\mbox{\rm on}\ \partial\Omega,$
\end{tabular} \right.
\endeq
where $\sigma\geq 0$ is a non-negative real number.

\begin{theorem}\label{texip}
Assume that $p,q\geq 0$ satisfy $p-q <1$.
\begin{itemize}
\item[\rm (i)] Then the system \eq{p} has solutions for all
$\sigma\geq 0$; \item[\rm (ii)] For any solution $(u,v)$ of
\eq{p}, there exist $c_1,c_2>0$ such that
\neweq{dit}
c_1 d(x)\leq u,v \leq c_2 d(x) \quad\mbox{ in } \Omega.
\endeq
Moreover, the following properties hold.
\begin{itemize}
\item[\rm (ii1)] If $-1<p-q<0$ then $u,v\in C^2(\Omega)\cap
C^{1,1+p-q}(\overline\Omega)$; \item[\rm (ii2)] If $0\leq p-q<1$
then $u,v\in C^2(\overline\Omega)$.
\end{itemize}
\end{itemize}
\end{theorem}
\noindent{\it Proof.} Existence follows directly from Theorem
\ref{texi} since conditions $(A1)$ and $(A2)$ are fulfilled.

\noindent (ii) Recall that from \eq{mzeta} we have $u\geq
\zeta\geq \overline c\varphi_1$ in $\Omega$. From the second
equation in \eq{p} we deduce
$$\di \Delta v-\beta v+\overline c^{p+\sigma}\frac{\varphi_1^{p+\sigma}}{v^{q+\sigma}}\leq 0
\quad \mbox{ in } \Omega.$$ Since $p-q<1$, we also get that
$\underline v= \underline c \varphi_1$ satisfies
$$\di \Delta \underline v-\beta \underline v+\overline c^{p+\sigma}\frac{\varphi_1^{p+\sigma}}{\underline v^{q+\sigma}}\leq 0
\quad \mbox{ in } \Omega,$$ provided $\underline c>0$ is
sufficiently small. Therefore, by virtue of Lemma \ref{maxp}, we
obtain $v\geq \underline c\varphi_1$ in $\Omega$.

Let us prove now the second inequality in \eq{dit}. To this aim,
set $w=u-v$. With the same idea as in Lemma \ref{lemep1} one gets
$\di \Delta w-\alpha w +\rho(x)\geq 0$ in the set $\{x\in\Omega:
w(x)>0\}$. Hence
\neweq{dubluw}
w\leq \zeta\leq c\varphi_1\quad \mbox{ in } \Omega.
\endeq
Let $w^+=\max\{w,0\}$. Then $v$ satisfies
$$
\left\{\begin{tabular}{ll} $\di \Delta v-\beta
v+\frac{(w^++v)^{p+\sigma}}{v^{q+\sigma}}\geq 0$
\quad & $\mbox{ in }\ \Omega,$\\
$\di v=0$\quad & $\mbox{ on }\ \partial\Omega.$\\
\end{tabular} \right.
$$
Consider now the problem
\neweq{duw}
\left\{\begin{tabular}{ll} $\di \Delta z-\beta
z+2^{p+\sigma}z^{p-q}=0$
\quad & $\mbox{ in }\ \Omega,$\\
$\di z>0$\quad & $\mbox{ in }\ \Omega,$\\
$\di z=0$\quad & $\mbox{ on }\ \partial\Omega.$\\
\end{tabular} \right.
\endeq
The existence of a classical solution to \eq{duw} follows from
\cite[Lemma 2.4]{sy1}. Moreover, if $0\leq p-q<1$ then $z\in
C^2(\overline\Omega)$ and with the same arguments as in
\cite[Theorem 1.1]{gl} we have $z\in C^2(\Omega)\cap
C^{1,1+p-q}(\overline\Omega)$ in the case $-1<p-q<0$.  Furthermore
$z\leq m\varphi_1$ in $\Omega$ for some $m>0$. On the other hand,
$\tilde c\varphi_1$ is a subsolution of \eq{duw} provided $\tilde
c>0$ is small enough. Therefore, by Lemma \ref{l1p} we get $z\geq
\tilde c\varphi_1$ in $\Omega$. This last inequality together with
\eq{dubluw}  allows us to choose $M>1$ large enough such that
$Mz\geq w^+$ in $\Omega$. Hence
$$\begin{aligned}
\di \Delta (Mz)-\beta
(Mz)+\frac{(w^++Mz)^{p+\sigma}}{(Mz)^{q+\sigma}} &
\leq \Delta (Mz)-\beta (Mz)+2^{p+\sigma}(Mz)^{p-q}\\
& \di = M\Big(\Delta z-\beta z+2^{p+\sigma}z^{p-q}\Big)=0 \quad
\mbox{\rm in }  \Omega.
\end{aligned}
$$
This means that $\overline v:=Mz$ verifies
$$\di \Delta \overline v-\beta \overline v+\frac{(w^+
+\overline v)^{p+\sigma}}{\overline v^{q+\sigma}}\leq 0 \quad
\mbox{\rm in } \Omega  \quad\mbox{\rm and }\quad \overline v=0
\quad \mbox{\rm on } \partial\Omega.
$$
Remark now that $\Psi(x,t)=-\beta
t+\frac{(w^+(x)+t)^{p+\sigma}}{t^{q+\sigma}}$,
$(x,t)\in\overline{\Omega}\times(0,\infty)$ satisfies the
hypotheses in Lemma \ref{l1p} since $p-q<1$. Furthermore, we have
$$\Delta \overline v+\Psi(x,\overline v)\leq 0\leq \Delta v+\Psi(x,v) \quad\mbox{\rm in }
\Omega,$$
$$\overline v,v>0  \quad\mbox{\rm in } \Omega, \quad
\overline v=v=0 \quad \mbox{\rm on }\partial\Omega \quad\mbox{\rm
and } \Delta \overline v\in L^1(\Omega).$$ Hence, by Lemma
\ref{l1p} we obtain
\neweq{dubluww}
\di v\leq \overline v\leq \tilde c\varphi_1 \quad \mbox{ in }
\Omega.
\endeq
Combining \eq{dubluw} and \eq{dubluww} we deduce $u=w+v \leq
C\varphi_1$ in $\Omega$, for some $C>0$. This completes the proof
of (ii1). As a consequence, there exists $M>1$ such that
$$\di 0 \leq \frac{u^{p}}{v^q},   \frac{u^{p+\sigma}}{v^{q+\sigma}}\leq M
\varphi_1^{p-q} \quad \mbox{ in } \Omega.$$ If $0\leq p-q<1$ then
by classical regularity arguments we have $u,v\in
C^2(\overline\Omega)$. If $-1<p-q<0$, then the same method as in
\cite[Theorem 1.1]{gl} can be employed in order to obtain $u,v\in
C^2(\Omega)\cap C^{1,1+p-q}(\overline\Omega)$.

This finishes the proof of Theorem \ref{texip}. \qed

\section{Uniqueness of the solution in one dimension}

In this section we are concerned with the uniqueness of the
solution associated to the one dimensional system
\neweq{1}
\left\{\begin{tabular}{ll} $\di u''-\alpha
u+\frac{u^p}{v^q}+\rho(x)=0$
\quad & $\mbox{\rm in }\ (0,1),$\\
$\di v''-\beta v+\frac{u^{p+\sigma}}{v^{q+\sigma}}=0$
\quad & $\mbox{\rm in }\ (0,1),$\\
$u(0)=u(1)=0, \; v(0)=v(1)=0.$
\end{tabular} \right.
\endeq
Our approach is inspired by the methods developed in
\cite{choimck1}, where a $C^2-$regularity of the solution up to
the boundary is needed. So, we restrict our attention to the case
$0<q\leq p\leq 1$. Thus, by virtue of Theorem \ref{texip}, any
solution of \eq{1} belongs to $C^2[0,1]\times C^2[0,1]$. By Hopf's
maximum principle we also have that $u'(0)>0$, $v'(0)>0$,
$u'(1)<0$ and $v'(1)<0$ for any solution $(u,v)$ of system \eq{1}.

The main result of this section is the following
\begin{theorem}\label{texiu}
Assume that $0<q\leq p \leq 1$, $\sigma\geq 0$. Then the system
\eq{1} has a unique solution $(u,v)\in C^2(\overline\Omega)\times
C^2(\overline\Omega)$.
\end{theorem}
\noindent{\it Proof.} Existence follows from Theorem \ref{texip}.
We prove here only the uniqueness. Suppose that there exist
$(u_1,v_1),(u_2,v_2) \in C^2[0,1]\times C^2[0,1]$ two distinct
solutions of \eq{1}.

First we claim that we can not have $u_2\geq u_1$ or $v_2\geq v_1$
in [0,1]. Indeed, let us assume that $u_2\geq u_1$ in [0,1]. Then
$$\di v''_2-\beta v_2+\frac{u_2^{p+\sigma}}{v_2^{q+\sigma}}=0=
v''_1-\beta v_1+\frac{u_1^{p+\sigma}}{v_1^{q+\sigma}} \quad
\mbox{\rm in }\ (0,1),$$ and by Lemma \ref{maxp} we get $v_2\geq
v_1$ in $[0,1]$. On the other hand
\neweq{doinw}
\di u''_1-\alpha u_1+\frac{u_1^{p}}{v_2^{q}}+\rho(x)\leq 0=
u''_2-\alpha u_2+\frac{u_2^{p}}{v_2^{q}}+\rho(x) \quad \mbox{\rm
in }\ (0,1).
\endeq
Note that the mapping $\Psi(x,t)=-\alpha
t+\frac{t^{p}}{v_2(x)^{q}}+\rho(x)$, $(x,t)\in
(0,1)\times(0,\infty)$ satisfies the hypotheses in Lemma \ref{l1p}
since $p\leq 1$. Hence $u_2\leq u_1$ in $[0,1]$, that is
$u_1\equiv u_2$. This also implies $v_1\equiv v_2$, contradiction.
Replacing $u_1$ by $u_2$ and $v_1$ by $v_2$, we also get that the
situation $u_1\geq u_2$ or $v_1\geq v_2$ in [0,1] is not possible.

\medskip

Set  $U=u_2-u_1$ and $V=v_2-v_1$. From the above arguments, both
$U$ and $V$ change sign in $(0,1)$. The key result in the approach
is the following.
\begin{prop}\label{psign}
$U$ and $V$ vanish only at finitely many points in the interval
$[0,1]$.
\end{prop}
\noindent{\it Proof.} We write the system \eq{1} as
$$
\left\{\begin{tabular}{ll} $\di {\bf W}''(x)+A(x){\bf W}(x)={\bf
0}$
\quad & $\mbox{\rm in }\ (0,1),$\\
${\bf W}(0)={\bf W}(1)={\bf 0},$
\end{tabular} \right.
$$
where ${\bf W}=(U,V)$ and $A(x)=(A_{ij}(x))_{1\leq i,j\leq 2}$ is
a $2\times 2$ matrix defined as
$$
A_{11}(x)=-\alpha+\left\{\begin{tabular}{cl} $\di
\frac{1}{v_2^q(x)}\cdot \frac{u_2^p(x)-u_1^p(x)}{u_2(x)-u_1(x)}\
,$
\quad & $ u_1(x)\neq u_2(x)$\\
$\di p \frac{u_1^{p-1}(x)}{v_1^q(x)}\ ,$
\quad & $ u_1(x)= u_2(x)$\\
\end{tabular} \right.
$$
$$
A_{12}(x)=\left\{\begin{tabular}{cl} $\di
-\frac{u_1^p(x)}{v_1^q(x) v_2^q(x)}\cdot
\frac{v_2^q(x)-v_1^q(x)}{v_2(x)-v_1(x)}\ ,$
\quad & $ v_1(x)\neq v_2(x)$\\
$\di -q \frac{u_1^{p}(x)}{v_1^{q+1}(x)}\ ,$
\quad & $ v_1(x)= v_2(x)$\\
\end{tabular} \right.
$$
$$
A_{21}(x)=\left\{\begin{tabular}{cl} $\di
\frac{1}{v_2^{q+\sigma}(x)}\cdot
\frac{u_2^{p+\sigma}(x)-u_1^{p+\sigma}(x)}{u_2(x)-u_1(x)}\ ,$
\quad & $ u_1(x)\neq u_2(x)$\\
$\di (p+\sigma) \frac{u_1^{p+\sigma-1}(x)}{v_1^{q+\sigma}(x)}\ ,$
\quad & $ u_1(x)= u_2(x)$\\
\end{tabular} \right.
$$
$$
A_{22}(x)=-\beta-\left\{\begin{tabular}{cl} $\di
\frac{u_1^{p+\sigma}(x)}{v_1^{q+\sigma}(x) v_2^{q+\sigma}(x)}\cdot
\frac{v_2^{q+\sigma}(x)-v_1^{q+\sigma}(x)}{v_2(x)-v_1(x)}\ ,$
\quad & $ v_1(x)\neq v_2(x)$\\
$\di (q+\sigma) \frac{u_1^{p+\sigma}(x)}{v_1^{q+\sigma+1}(x)}\ ,$
\quad & $ v_1(x)= v_2(x)$\\
\end{tabular} \right.
$$
Therefore, $A \in C(0,1)$ and  $A_{12}(x)\neq 0$, $A_{21}(x)\neq
0$ for all $x\in (0,1)$. Moreover, $xA(x)$, $(1-x)A(x)$ are
bounded in $L^\infty(0,1)$. Indeed, let us notice first that, by
\eq{dit} in Theorem \ref{texip}, there exist $c_1,c_2>0$ such that
$$\di c_1\leq \frac{u_i}{\min \{x,1-x\}}, \frac{v_i}{\min \{x,1-x\}} \leq c_2,\; (i=1,2)
\quad \mbox{\rm in }\ (0,1).$$ Then, by the mean value theorem, we
have
$$
\begin{aligned}
\di x |A_{12}(x)|& \leq qx \frac{u_1^p(x)}{v_1^q(x) v_2^q(x)}
\max \{v_1^{q-1}(x), v_2^{q-1}(x)\}\\
\di & \leq qx^{p-q} \left(\frac{u_1(x)}{x}\right)^p \max \left\{
\left(\frac{x}{v_1(x)}\right)^{q+1},
\left(\frac{x}{v_2(x)}\right)^{q+1} \right\}\\
\di & \leq c x^{p-q} \quad \mbox{ for all }\ 0<x\leq 1/2.
\end{aligned}
$$
We obtain similar estimates for $x A_{11}$, $x A_{21}$ and $x
A_{22}$. This allows us to employ Lemma 7 and Lemma 8 in
\cite{choimck1}. Note that condition $x A(x)\in L^\infty(0,1)$
suffices in order to obtain the same conclusion as in \cite[Lemma
8]{choimck1}. In particular, we get that $U$ and $V$ vanish only
at finitely many points in any compact interval
$[a,b]\subset(0,1)$.

It remains to show that $U$ and $V$ can not have infinitely many
zeroes in the neighborhood of $x=0$ and $x=1$. We shall consider
only the case $x=0$, the situation where $u$ or $v$ vanishes for
infinitely many times near $x=1$ being similar.

Without loosing the generality, we may assume that $U$ has
infinitely many zeroes in a neighborhood of $x=0$. Since $U\in
C^2[0,1]$ by Rolle's Theorem we get that both $U'$ and $U''$ have
infinitely many zeros near $x=0$. As a consequence, we obtain
$U'(0)=0$, that is, $u'_1(0)=u'_2(0)$.

If $V'(0)=0$, then ${\bf W}(0)={\bf W}'(0)={\bf 0}$ and by
\cite[Lemma 8]{choimck1} we deduce ${\bf W}\equiv {\bf 0}$ in
$[0,1/2]$ which is a contradiction. Hence $V'(0)\neq 0$.
Subtracting the first equation in \eq{1} corresponding to $u_1$
and $u_2$ we have
$$
\begin{aligned}
\di U''(x)& \di =\alpha U(x)+\frac{u_1^p(x)}{v_1^q(x)}-\frac{u_2^p(x)}{v_2^q(x)}\\
&\di =x^{p-q}\left\{\alpha
\frac{U(x)}{x^{p-q}}+\left(\frac{u_1(x)}{x}\right)^p
\left(\frac{x}{v_1(x)}\right)^{q}-
\left(\frac{u_2(x)}{x}\right)^p\left(\frac{x}{v_2(x)}\right)^{q}\right\}.
\end{aligned}
$$
Since $0\leq p-q<1$, $u'_1(0)=u'_2(0)$ and $v'_1(0)\neq v'_2(0)$
we get
$$\lim_{x\ri 0^+}\left\{\alpha \frac{U(x)}{x^{p-q}}+\left(\frac{u_1(x)}{x}\right)^p
\left(\frac{x}{v_1(x)}\right)^{q}-
\left(\frac{u_2(x)}{x}\right)^p\left(\frac{x}{v_2(x)}\right)^{q}\right\}=
u'^p_1(0)\left(\frac{1}{v'^q_1(0)}-\frac{1}{v'^q_1(0)}\right)\neq
0.$$ Therefore, $U''$ has constant sign in a small neighborhood of
$x=0$ which contradicts the above arguments. The proof of
Proposition \ref{psign} is now complete. \qed
\medskip

{\bf Proof of Theorem \ref{texiu} completed}.

Set
$$\di {\mathcal I}^+=\{x\in [0,1]: U(x)\geq 0\},\quad
\di {\mathcal I}^-=\{x\in [0,1]: U(x)\leq 0\},$$
$$\di {\mathcal J}^+=\{x\in [0,1]: V(x)\geq 0\},\quad
\di {\mathcal J}^-=\{x\in [0,1]: V(x)\leq 0\}.$$ According to
Proposition \ref{psign}, the above sets consist of finitely many
disjoint closed intervals. Therefore, ${\mathcal
I}^+=\cup_{i=1}^{m} I^+_i$. For simplicity, let $I^+$ denote any
interval $I^+_i$ and we use similar notations for $I^-$, $J^+$ and
$J^-$. We have
\begin{lem}\label{lintv}
For any intervals $I^+$, $I^-$, $J^+$ and $J^-$ defined above, the
following situations can not occur:

${\rm (i)}\;\;  I^+\subset J^+$; $\;\quad{\rm (ii)}\;\; I^-\subset
J^-$; $\;\quad{\rm (iii)}\;\;J^+\subset I^-$; $\;\quad {\rm
(iv)}\;\;J^-\subset I^+.$
\end{lem}
\noindent{\it Proof.} (i) Assume that $I^+\subset J^+$. Since
$v_2\geq v_1$ in $I^+$ we deduce that the inequality \eq{doinw}
holds in $I^+$. Using the fact that $u_2=u_1$ on $\partial I^+$,
by virtue of Lemma \ref{l1p} we get $u_2\leq u_1$ in $I^+$. Hence,
$u_2\equiv u_1$ in $I^+$, which contradicts Proposition
\ref{psign}. Similarly we can prove the statement (ii).

(iii) Assume that $J^+\subset I^-$. Then
$u_1^{p+\sigma}/v_1^{q+\sigma}\geq u_2^{p+\sigma}/v_2^{q+\sigma}$
in $J^+$. Notice that $V=v_2-v_1$ verifies
$$
\left\{\begin{tabular}{ll} $\di V''-\beta
V=\frac{u_1^{p+\sigma}}{v_1^{q+\sigma}}-
\frac{u_2^{p+\sigma}}{v_2^{q+\sigma}}\geq 0$
\quad & $\mbox{\rm in }\ J^+,$\\
$\di V=0$\quad & $\mbox{\rm on } \partial J^+.$
\end{tabular} \right.
$$
By the maximum principle, it follows that $V\leq 0$ in $J^+$,
i.e., $v_2\leq v_1$ in $J^+$. This  yields $v_2\equiv v_1$ in
$J^+$ which again contradicts Proposition \ref{psign}. The proof
of (iv) follows in the same manner. \qed
\medskip

From now on, the proof of Theorem \ref{texiu} is the same as in
\cite[Theorem 6]{choimck1}. \qed

\medskip

\begin{figure}[tbph]
\begin{center}
\includegraphics[height=6cm,angle=0]{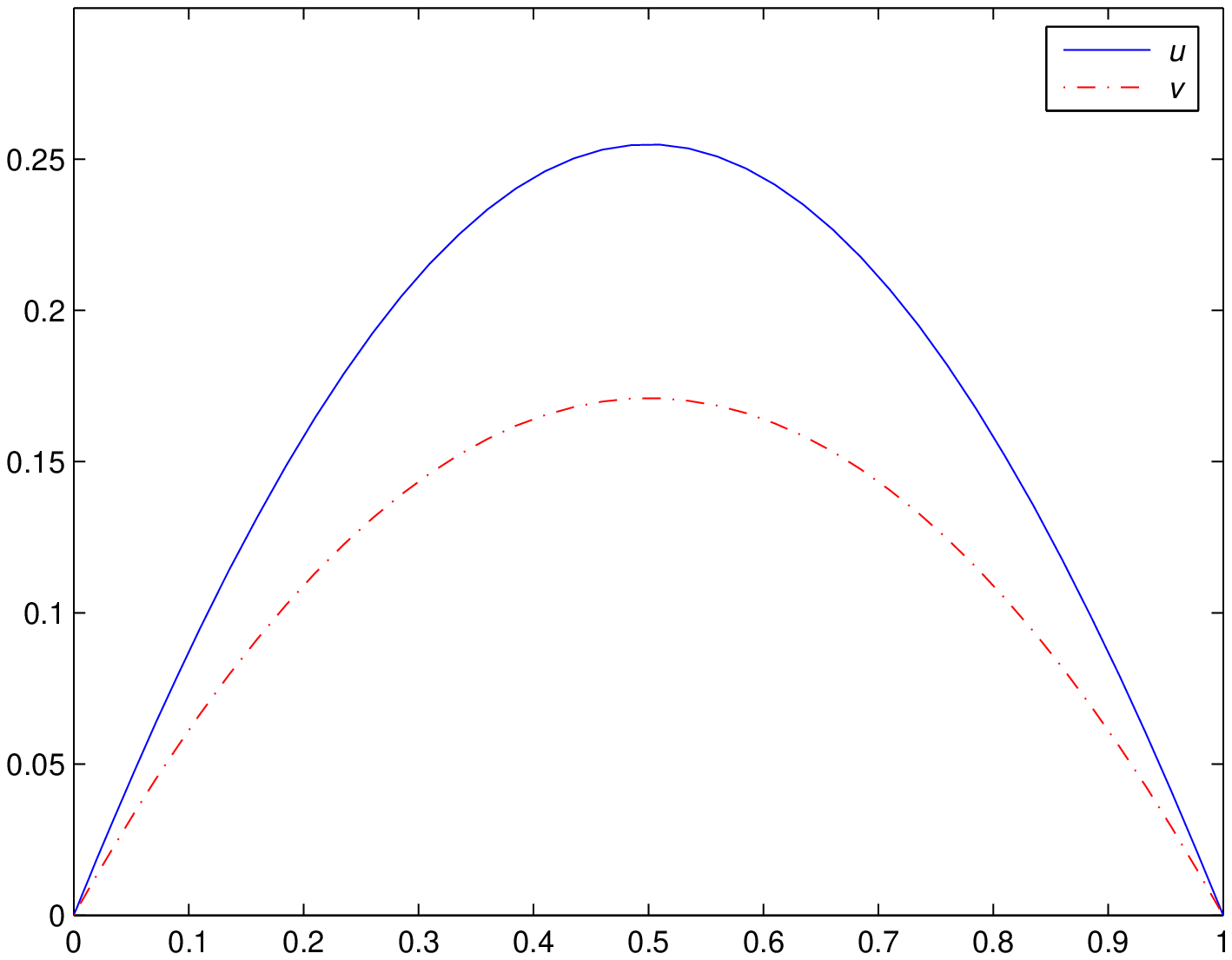}
\includegraphics[height=6cm,angle=0]{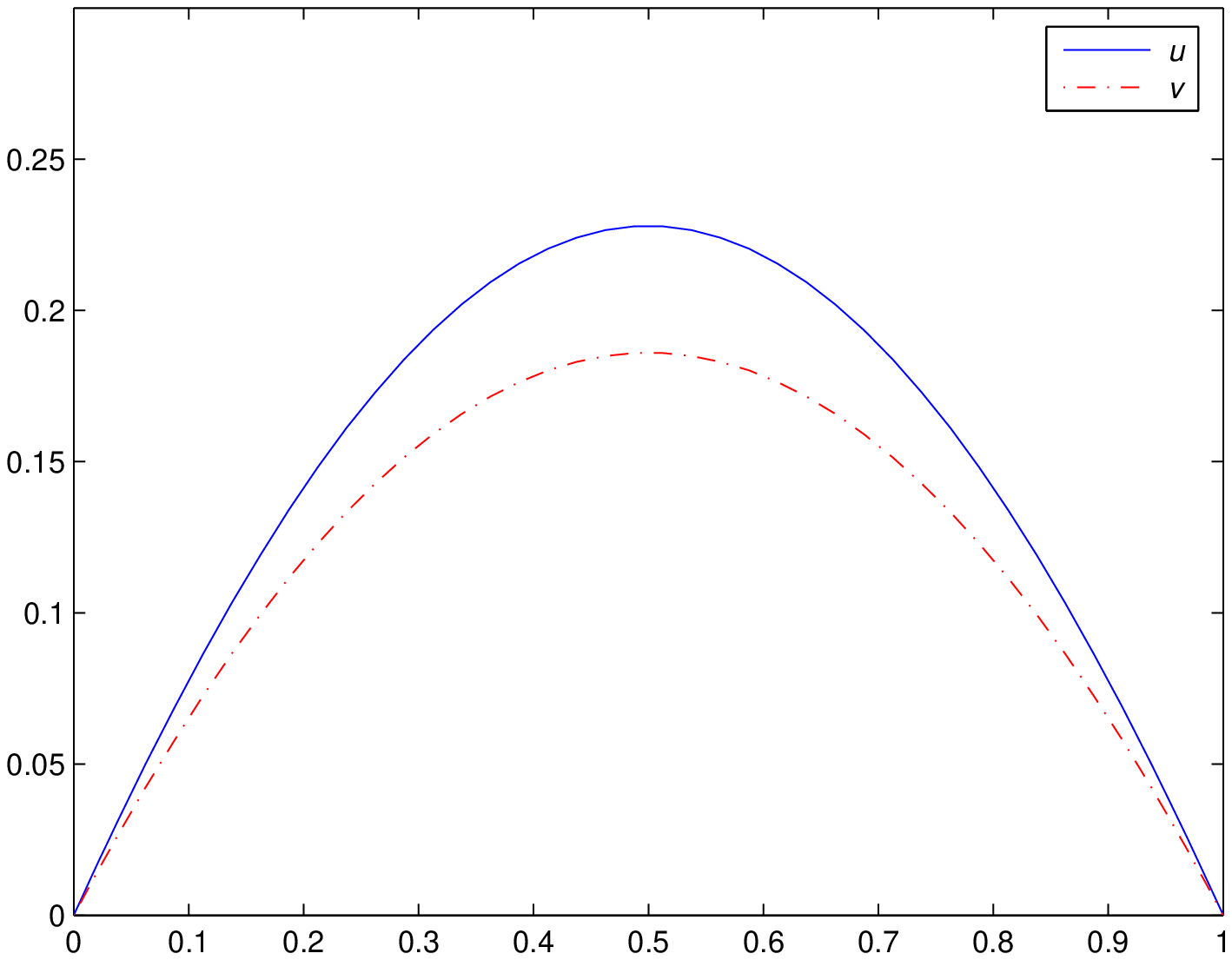}
\caption{The solution $(u,v)$ of the system $(S_\ep)$ with
$\alpha=1$, $\beta=0.5$, $p=q=1$, $\ep=10^{-2}$ and
$\rho(x)=\sin(\pi x)$. We have chosen $\sigma=0$ (on the left) and
$\sigma=2$ (on the right).}
\end{center}\label{buu}
\end{figure}
\medskip

\noindent{\bf Remark.} As a consequence of Theorem \ref{texi}, the
solution $(u,v)$ of the system \eq{p} can be approximated by the
solutions of $({\mathcal S})_\ep$. Furthermore, the shooting
method combined with the Broyden method in order to avoid the
derivatives, are suitable to numerically approximate the solution
of \eq{p}. We have considered $\alpha=1$, $\beta=0.5$, $p=q=1$,
$\ep=10^{-2}$ and $\rho(x)=\varphi_1(x)=\sin(\pi x)$. In the above
figure we have plotted the solution $(u,v)$ of $(S_\ep)$ for
$\sigma=0$ (on the left) and $\sigma=2$ (on the right)
respectively.

\noindent{\bf Acknowledgement.} The authors are partially
supported by Grant 2-CEx06-11-18/2006.

\end{document}